\documentclass[12pt,a4paper]{article}
\usepackage{amsmath, amssymb, theorem, latexsym,color}
\usepackage{graphicx}

\newtheorem{theorem}{Theorem}[section]

\newtheorem{proposition}[theorem]{Proposition}
\newtheorem{definition}[theorem]{Definition}

\newtheorem{corollary}[theorem]{Corollary}

\newtheorem{remark}[theorem]{Remark}

\allowdisplaybreaks

\def\neweq#1{\begin{equation}\label{#1}}
\def\endeq{\end{equation}}
\newcommand\finbox{~\hfill$\Box$}%

\newcommand {\ar}{\rightarrow}

\numberwithin{equation}{section}

\title{
Pleijel's theorem for Schr\"odinger operators with radial potentials}

\author{Philippe Charron\\ 
Universit\'e de Montr\'eal\,\\
2920, Chemin de la Tour, 
Montr\'eal, QC,
H3T 1J4, Canada\\~\\
Bernard Helffer\\
 Laboratoire de Math\'ematiques Jean Leray, Universit\'e de Nantes\\
 2 rue de la Houssini\`ere 44322 Nantes, France\\
 and\\
 Laboratoire de Math\'ematiques, \\ Universit\'e Paris-Sud, CNRS, Univ. Paris Saclay, France\\~\\
and\\~\\
 Thomas Hoffmann-Ostenhof\\
 Vienna University, Department of Theoretical Chemistry,\\
A 1090 Wien, W\"ahringerstrasse 17, Austria}
 
\date{}
\begin{document}
\maketitle
\begin{abstract}
In 1956, $\AA$. Pleijel gave his celebrated theorem showing that the inequality in Courant's theorem on the number of nodal domains is strict for large eigenvalues of the Laplacian. This was a consequence of a stronger result giving an asymptotic upper bound for the number of nodal domains 
 of the eigenfunction as the eigenvalue tends to $+\infty$.   A similar question occurs naturally for the case of the Schr\"odinger operator. The first  significant result has been  obtained recently by the first author for the case of the harmonic oscilllator. The purpose of this paper is to consider more general potentials which are radial. We will analyze either the case when the potential tends to $+\infty$ or the case when the potential tends to zero, the considered eigenfunctions being associated with the eigenvalues below the  essential spectrum.
\end{abstract}

\section{Introduction}
The goal of this paper is to extend Pleijel's theorem for the Dirichlet Laplacian $H(\Omega)=-\Delta$ in a bounded domain $\Omega$  to the case of the Schr\"odinger operator $H_V=-\Delta +V$ in $\mathbb R^d$.  We are interested in counting the number of nodal domains of an eigenfunction and to relate this number with the labelling of the corresponding eigenvalue. Throughout this paper, for any function $f$ defined over a domain $D\subset \mathbb R^d$, $\mu(f)$ will denote the number of nodal domains of $f$, namely the number of connected components of $D \setminus f^{-1}(0)$. The starting point of the analysis is Courant's Theorem (1923)  \cite{Cou}.
 \begin{theorem}[Courant]
If $\phi_n$ is an eigenfunction associated with the $n$-th eigenvalue  $\lambda_n$  of $H(\Omega)$ (ordered in non decreasing order and labelled with multiplicity), then
\begin{equation}\label{CourIn}
\mu(\phi_n) \leq n\,.
\end{equation}
\end{theorem}
Pleijel's theorem (1956) \cite{Pl}  says
\begin{theorem}[Pleijel's weak theorem]\label{Pleijelweakform}
If the dimension is $\geq 2$, there is only a finite number of eigenvalues of the Dirichlet Laplacian in $\Omega$ for which we have equality in \eqref{CourIn}.
\end{theorem} 
Let us now give the strong form of Pleijel's theorem.  
\begin{theorem}[Pleijel's strong theorem]\label{PropPleijel1}
Let $\lambda_n$ the non decreasing sequence of eigenvalues associated to the Dirichlet realization of the Laplacian.
For any $d \geq 2$   for any orthonormal
basis $\phi_{n}$ of eigenfunctions $\phi_n$  of $H(\Omega)$ (the Dirichlet Laplacian in $\Omega$) associated with $\lambda_{n}(\Omega)$,
 \begin{equation}\label{bourgain3} 
\limsup_{n\ar +\infty} \frac{\mu(\phi_n)}{n}\leq \gamma (d)  =  \frac{2^{d-2}d^2 \Gamma(d/2)^2}{(j_{\frac{d}{2}-1})^d},
	\end{equation}
	where $ j_\nu$ denotes the first zero of the Bessel function $ J_\nu$.
\end{theorem}
 The theorem was proved by Pleijel \cite{Pl}  for $d=2$ and then extended by Peetre \cite{Pe}  and B\'erard-Meyer \cite{BeMe}. 
We recall from \cite[Lemma~9]{BeMe} that  Pleijel's constant equals
\begin{equation}\label{gamma2}
\gamma (d)= w_{d}^{-1} \omega_d^{-1} \, \bigl(\lambda(B_d) \bigr)^{- d/2} < 1\,,
\end{equation}
where
\begin{itemize}
\item  $w_{d}$ is the Weyl constant 
\begin{equation}\label{defw}
w_{d}:= (2\pi)^{-d}\omega_d\,,
\end{equation}
\item 
\begin{equation}\label{defomega}
\omega_d:= |B_d|\,,
\end{equation}
where $B_d$ is  the unit ball in $\mathbb R^d$ and   $|D|$ denotes for an open set  $D\subset \mathbb R^d$  its volume;
\item $\lambda(B_d) $ is the Dirichlet ground state energy of the 
Laplacian in $B_d$.
\end{itemize}

As  $\gamma(d) <1$,  one recovers as a corollary that the inequality in Courant's theorem is strict for $n$ large.   The second point to notice is that the constant is independent of the open set.  Complementary properties of $\gamma (d)$ have been obtained by B. Helffer and M. Persson-Sundqvist \cite{HPS}. In particular $d \mapsto \gamma (d)$ is decreasing exponentially  to $0$.  Finally note that this constant is not optimal (see \cite{Bo}, \cite{St} and the discussion in \cite{HHO}).

The original proof of Pleijel's theorem is based on a combination between the Weyl formula and the Faber-Krahn inequality.
Weyl's theorem  reads, as $\lambda \ar +\infty$, 
\begin{equation}\label{WF}
N (\lambda)  = |\Omega|\, w_d  \,  \lambda^{\frac d2} (1 + o(1))\,,
\end{equation}
where 
\begin{equation}\label{Nlambda}
N (\lambda) := \# \{ \lambda_j <\lambda\}\,.
\end{equation}
 
 The Faber-Krahn inequality  reads
 \begin{theorem}\label{ThFK}
For any domain $D\subset \mathbb R^d$ ($d\geq 2$), we have 
\begin{equation}\label{eq.FK}
|D|^\frac 2 d \,  \lambda(D) \geq \omega_d^{\frac 2d} \lambda(  B_d )\,.
\end{equation}
\end{theorem}

There are a lot of Weyl's formulas available in the context  of the Schr\"odinger operator $ H_V:=-\Delta + V\,.$
 The use of the Faber-Krahn inequality is more problematic, except  of course for  the case of bounded domains with
  bounded potential which can be treated like the membrane case. In 1989 Leydold \cite{Ley}  obtained in his diploma thesis a
  weak Pleijel theorem for the isotropic harmonic oscillator (see also \cite{BeHe}).
  Two years ago Charron \cite{Cha} in his master thesis proved Pleijel's strong
  theorem also for the
  harmonic oscillator:
\begin{theorem}[Charron's theorem]
	Let $(\phi_{n})_{n\in \mathbb N}$ be an orthogonal basis of eigenfunctions in $L^2(\mathbb R^d)$  of the harmonic oscillator associated with $(\lambda_{n})_{n\in \mathbb N}\,$. Then 
	\begin{equation}\label{bourgain3} 
	\limsup_{n\ar +\infty} \frac{\mu(\phi_n)}{n}\leq\gamma (d)\,.
	\end{equation}
	\end{theorem}

The theorem is also proven in the case of the non-isotropic harmonic oscillator \cite{Cha1}.  The interesting fact is that the potential $V(x) =\sum_{i} a_{i}x_i^2$ (with $a_i >0$) does not appear  on the right hand side of the upper bound. Note also that when there are no eigenvalue degeneracies a much stronger result is available in \cite{Cha1}.\\

A natural question to ask is whether the theorem can be extended to more general Schr\"odinger operators. We will answer positively this question under the additional assumption that the potential is radial.\\
More precisely, we assume
$$d \geq 2$$
and  we consider on $\mathbb{R}^d$ a Schr\"odinger operator 
$
{H_V} = - \Delta + V \, ,$
where $V$ is  a radial potential:
\begin{equation}\label{hyprad}
V(x) = v(r)\,,
\end{equation}
with $|x|=r$.\\
We will assume that 
\begin{equation}\label{ass1}
v\in C^1 (0,+\infty)\,,
\end{equation}
and that there exists $R_0> 0$ such that
\begin{equation}\label{assmonotone}
v'(r) >0\,, \mbox{ for } r\geq R_0\,.
\end{equation}
In order to allow some singularity at the origin, 
we assume
either 
\begin{equation}\label{ass3-0}
v \in C^0([0,+\infty))\,,
\end{equation}
or that there exists $s \in (0,2)$ such that, as $r\ar 0\,$, 
\begin{equation}\label{ass3}
v(r) \approx -  r^{-s} \,.
\end{equation}
Here, we say that $a\approx b$ is $a/b$ and $b/a$ are bounded, i.e. if there exists $C >0$ such that
$$
\frac 1C \leq \frac a b \leq C\,.
$$
 We will study two cases according to the behavior of $v$ at $+\infty\,$.\\
 {\bf Case A:} 
$v$  tends to $+\infty$  as $r \to + \infty\,$. \\
More precisely, we assume  \eqref{ass1}-\eqref{ass3}  and that there exists\footnote{This condition appears when applying Weyl's 
formula given by Theorem \ref{WRS} from \cite{ReSi}. At least under stronger assumptions on the regularity of $v$ 
for  $r\rightarrow +\infty$, it should be possible to assume $m>0$.} $m >1$ such that as $r\ar +\infty$
\begin{equation}\label{ass4}
v(r) \approx r ^m\,,
\end{equation}
and
\begin{equation}\label{ass5}
v'(r) \approx r^{m-1}\,.
\end{equation}
{\bf Case B}:
$v $ tends to $0$  as $ r \to +\infty\,$. \\
More precisely, we assume  \eqref{ass1}-\eqref{ass3} and that there exists $m \in (-2,0)$ such that
\begin{equation}\label{ass6}
v(r) \approx - r ^m\,,
\end{equation}
and
\begin{equation}\label{ass7}
v'(r) \approx r^{m-1}\,.
\end{equation}

 In the two cases there is a natural selfadjoint extension starting from $C_0^\infty(\mathbb R^d)$ (see Section \ref{s3}). 
In  Case A, the spectrum is discrete and  consists of a non decreasing sequence of eigenvalues $\lambda_n$  tending to $+\infty\,$.  In  Case~B    the spectrum is divided in two parts, the essential spectrum:  $[0,+\infty)$ 
 and the discrete spectrum, which consists of an infinite sequence  of negative eigenvalues  $(\lambda_n)_{n\in \mathbb N}$ tending to $0$ (see for example Reed-Simon \cite{ReSi}, Vol. IV, Theorem~XIII.6). 
Associated with this sequence $(\lambda_n)_{n\in \mathbb N}$, we can consider an orthonormal sequence of eigenfunctions $\phi_n$, where in Case A    $\phi_n$ is an Hilbertian basis of $L^2(\mathbb R^d)$ and in Case B    of the negative eigenspace.

Our analysis will contain  two well-known potentials: the quantum harmonic oscillator (Case A) and the Coulomb potential (Case B). In both cases, we know the eigenvalues and  an explicit basis of eigenfunctions  but in the proof this property will not be used. Our aim is to prove the following result: 
 \begin{theorem}[Pleijel's theorem for Schr\"odinger]\label{thm1.6}
 In  Cases  A     or   B,  if  $(\phi_{n})_{n\geq 1}$ is an orthogonal  sequence of eigenfunctions of $H_V$ associated with the above defined sequence $\lambda_{n}$, then 
	\begin{equation}\label{bourgain3} 
	\limsup_{n\ar +\infty} \frac{\mu(\phi_n)}{n}\leq \gamma (d)\,.
	\end{equation}
\end{theorem}

The paper is organized as follows.\\
In Section \ref{s2}, we discuss the general strategy and the methods used by Pleijel first and then by P. Charron.
In Section  \ref{s3}, we review the general properties of the Schr\"odinger operator.
In Section \ref{s4}  we collect those Weyl-type 
results we need for
the proof of Theorem \ref{thm1.6}. 

In Section \ref{s5}, we give the proof of  our Pleijel's theorem in the two situations.
\\

 {\bf Acknowledgements.}\\
Thanks to the ESI where the paper was initiated (B. H. and T. H.-O.). The authors would also thank I. Polterovich for helpful discussions at various stages of this work.

 \section{About the methods}\label{s2}
As recalled in the introduction, the original proof of $\AA$. Pleijel was based on a tricky combination of Weyl's formula with the Faber-Krahn inequality.  When considering the case of the Schr\"odinger operator in $\mathbb R^d$, Weyl's formula still exists but  the use of Faber-Krahn is not easy: nodal domains could be unbounded and the variation of the potential inside a nodal domain could be very high. One has consequently to find an idea for proving that these two bad situations do not occur very often.

In the case of the harmonic oscillator Charron's proof relies on specific properties of the eigenfunctions and the potential. Namely, it used the fact that every eigenfunction is a linear combination of an exponential multiplied by polynomials whose  degree can be controlled by a  function of the labeling of the eigenvalue, that the hypersurfaces with constant potential are hyperspheres and the fact that the counting function $N (\lambda) $
 behaves nicely as $\lambda \ar +\infty$ (Weyl's law).

In addition,  it also used that, every nodal domain of an eigenfunction of a  Schr\"odinger operator intersects the classically allowed region associated with the eigenvalue $\lambda$, i.e 
\begin{equation}\label{defenergy}
V^{(-1)}(-\infty,\lambda):=\{ x\in \mathbb R^d\,|\, V (x) < \lambda\}\,.
\end{equation}
This property is quite general and elementary.

The key was then to divide the classically allowed region in a finite number of annuli of the form $ V^{(-1)}(a,b)$. Every nodal domain can either be contained in a single annulus or intersect more than one. To give an upper bound on the number of nodal domains not  contained in one annulus, Charron uses properties of  algebraic surfaces, as well as results arising in Morse theory adapted from Milnor \cite{Mi}.

Then, he used the  Faber-Krahn inequality to give a lower bound on the volume of any nodal domain contained in a single annulus. Dividing the volume of each annulus by the volume of each nodal domain gave an upper bound on the number of nodal domains contained in that annulus.

The last step was to find an appropriate number of annuli to balance out both estimates.\\

To extend the methods of Charron's proof to more general potentials, we need to find Schr\"odinger operators such that: 
\begin{enumerate}
	\item There are good lower bounds  for the number of eigenvalues below any $\lambda$.
	\item We can count the number of nodal domains that intersect a given energy hypersurface $V^{(-1)} (b)$. 
	\item{We can give an upper bound on the number of nodal domains that are not contained in the classically allowed region.}
\end{enumerate}

In the case of the harmonic oscillator, the eigenvalues are known explicitly. However, for many potentials $V$, Weyl's  law  can be extended to the Schr\"odinger operator $H_V$ for estimating the number of eigenvalues. Hence, we need to find a suitable class of potentials where this  law holds. 

So far, the only known method to give a suitable upper bound on the number of nodal domains that intersect an energy hypersurface  are based on Milnor's theorem (see Subsection \ref{ss3.7}).  Hence we need  this hypersurface to be  algebraic and  the property that for any eigenvalue $\lambda$  and any energy hypersurface (or at least a suitable family)   the restriction of any associated eigenfunction  $u_\lambda $ equals the restriction of a polynomial to this hypersurface. This is why we  focus in this paper on the study of radial potentials.  In this case, the energy hypersurfaces are hyperspheres $\{ r = \rho\}$ for some $\rho>0$  and it can be shown that the restriction of an eigenfunction to a hypersphere is always a linear combination of hyperspherical harmonics, each one being the restriction to the hypersphere of a homogeneous harmonic polynomial.  We will also  have to control the degree of such polynomials by a  function of $\lambda$ or of its labelling.
This  last property will allow us to bound the number of nodal domains that are not contained in $V^{(-1)}(-\infty,\lambda)$.

Another problem might arise  when estimating  the number of nodal domains contained in one annulus. In the case of the harmonic oscillator, summing over all annuli gives us an expression which can be compared directly with an integral. The error term that arises becomes negligible as $\lambda \to + \infty\,$. It remains to show under which conditions on $V$ the same method can be applied.

Finally, in the specific case of Coulomb-like potentials at the origin, we need to look at the behavior of the number of nodal domains near  the origin.

\section{On the spectral theory of the Schr\"odinger operators  with radial potential}\label{s3}
\subsection{General theory}
We first verify that our Schr\"odinger operator $H_V =-\Delta +V$  is well defined by a Friedrichs procedure  starting from its sesquilinear form
defined on \break $C_0^\infty(\mathbb R^d \setminus \{0\})\times C_0^\infty(\mathbb R^d \setminus \{0\}) $
$$
(u,v)  \mapsto a(u,v) :=\int_{\mathbb R^d} \nabla u (x) \cdot \nabla v (x)   \, dx + \int_{\mathbb R^d}  V(x) \, u(x)\,v(x)\, dx\,.
$$
Note that the left term has a meaning as soon as $V\in L^1_{loc} (\mathbb R^d\setminus \{0\})\,$.  In our case, this is a consequence of Assumption \ref{ass1}.
Our operator, will be defined through a Friedrichs extension.  This  works as soon as  $q(u):= a(u,u)$ is bounded from below by $- C \,||u||^2_{L^2}\,$. 
  It is consequently enough to control
 the integral $\int_{ V < 0} V(x) \, u(x)^2 \, dx$ from below.\\
 When $d\geq 3$, we use Hardy's inequality
\begin{equation} \label{Hardyineq}
 \int _{\mathbb R^d}|\nabla u(x)|^2 dx \geq \left(\frac{d-2}{2} \right)^2 \int _{\mathbb R^d} \frac{1}{r^2}\, |u(x)|^2 \,dx \,,\, \forall u \in C_0^\infty(\mathbb R^d)\,,
\end{equation}
with $r=|x|\,$.\\
This inequality extends to $H^1(\mathbb R^d)$ by density.\\
For $d=2$, we can use the modified Hardy inequality  in a disk $D(0,\check R)$ which reads (see  for example \cite{Cow})
\begin{equation}\label{logHardy}
\int_{\mathbb R^2} |\nabla u (x)|^2 dx \geq \frac{1}{4}\int_{\mathbb R^2} \frac{u(x)^2}{|x|^2\log ^2(|x|/\check R)}\,dx\,,\, \forall u \in C_0^\infty(D(0,\check R))\,,
\end{equation}
which also extends to $H^1_0(D(0,\check R))\,$.\\
Using these inequalities and a partition of unity, the semi-boundedness on $C_0^\infty(\mathbb R^d\setminus\{0\})$ follows immediately.\\

Let us now  describe the form domain resulting from the Friedrichs extension procedure. We have:
\begin{itemize}
\item[ \bf Case A ]  
$$ 
Q_H =  \{u \in H^1(\mathbb R^d) \,|\, \sqrt{V} u \in L^2( D(0,R_1)^c)\}
$$
where $R_1$ is chosen such that $v(r) \geq 1$ for $r\in (R_1,+\infty)\,$,\\
\item[\bf Case B  ]  
$$
Q_H = H^1(\mathbb R^d)\,.
$$
\end{itemize}
We do not need to characterize the domain of the corresponding self-adjoint operator.

 \subsection{Nodal domains intersect the classically allowed region}\label{ss3.2}
We use a similar argument as in  \cite{Ley} and \cite{Cha1}\,.  We assume that we are either in Case A or in Case B, but the result is much more general. 
\begin{proposition} \label{lemma1}
Let $\lambda$ be an eigenvalue below the essential spectrum,  $u_\lambda $ be an eigenfunction of $H_V$ associated with eigenvalue $\lambda$ and $\Omega$ be a nodal domain of $u_\lambda $. 
Then
$$\Omega \cap V ^{(-1)}(-\infty,\lambda) \neq \emptyset\,.
$$
\end{proposition}
{\bf Proof.}\\
If for all $x \in \Omega$, $V(x) > \lambda\,$, then
\begin{multline}
	\lambda \, = \frac{\int_\Omega {{|\nabla u_\lambda (x)|}^2\, dx } + \int_\Omega {V(x)u_\lambda(x)^2\,dx }}{\int_\Omega u_\lambda (x)^2 \, dx} \\
	\quad \geq \frac{\int_\Omega {V(x)u_\lambda(x)^2}\, dx }{\int_\Omega {u_\lambda ( x)^2}\,dx }\\
	\quad > \frac{\int_\Omega {\lambda u_\lambda  (x)^2\,dx}}{\int_\Omega {u_\lambda (x)^2}\,dx} = \lambda \, ,
	\end{multline}
hence a contradiction.
\finbox

Therefore, any nodal domain is either contained in the classically allowed region $\left\lbrace V < \lambda\right\rbrace $ or intersects the hypersurface $V^{(-1)}(\lambda)\,$.

\subsection{The radial Schr\"odinger operator}
Although the exposition there is limited to the case $d=3$, one can refer to Reed-Simon \cite{ReSi}, Vol. IV p. 90-91.\\
The Laplace operator $-\Delta$ can be written as
\begin{equation}
-\Delta=-\frac{\partial^2}{\partial r^2}-\frac{d-1}{r}\frac{\partial}{\partial r}
+\frac{1}{r^2}(-\Delta_{\mathbb S^{d-1}}),
\end{equation}
where $r=|x|$ is the radial variable and $\Delta_{\mathbb S^{d-1}}$ is the
Laplace--Beltrami operator, acting in $L^2(\mathbb S^{d-1})$. The following proposition is standard (see for example \cite{Sh},
Theorem~22.1 and Corollary~22.1).

\begin{proposition}\label{proshu} 
Assume that $d\geq 2$. 
The spectrum of $-\Delta_{\mathbb S^{d-1}}$ consists of eigenvalues
\[
\ell(\ell+d-2)\,,\quad \ell\in\mathbb{N}\,.
\]
The multiplicity of the eigenvalue $\ell(\ell+d-2)$ is given by
\[
\Lambda_{\ell,d}:=\binom{\ell+d-1}{d-1}-\binom{\ell+d-3}{d-1}\,,
\]
which coincides with the dimension of the space of homogeneous, 
harmonic polynomials of degree $\ell$.
\end{proposition}
We denote by $\mathbb S^{d-1}\ni \omega \mapsto Y_{\ell,m}(\omega) $ an orthonormal basis of the $\Lambda_{\ell,d}$-dimensional eigenspace associated with  $\ell(\ell+d-2)$. We recall that each  $Y_{\ell,m}$ is the restriction to $\mathbb S^{d-1}$ of a harmonic homogeneous  polynomial of degree $\ell$.\\
We now consider the Schr\"odinger operator $H_V$ and assume
\begin{equation}
V(x) =v(r)\,.
\end{equation}
 In this case, one can determine the spectrum by  using polar coordinates. 
 In the spherical coordinates, we can determine the spectrum by considering the (closure of the) union of the spectra of
 the family (indexed by $\ell \in \mathbb N$) 
of Sturm Liouville 
operators $\mathcal L_\ell$ defined by
\begin{equation}\label{Ll}
\mathcal L_\ell=-\frac{d^2}{dr^2}-\frac{d-1}{r}\frac{d}{dr}+\frac{\ell(\ell+d-2)}{r^2} + v(r)\,,
\end{equation}
acting in $L^2((0,+\infty),r^{d-1}\,dr)\,$, with a suitable Dirichlet like condition at $0$ (see Reed-Simon \cite{ReSi}, p. 91, Proof of Lemma 1). Note that the "Dirichlet like" condition is expressed after the unitary transform 
 $u \mapsto r^{\frac{d-1}{2}} u$ sending \break $L^2((0,+\infty); r^{d-1} dr)$ onto $L^2((0,+\infty); dr)$ and  becomes the standard Dirichlet condition for $\ell=0$. When $\ell >0$, no condition is given. The new operator is then:
 \begin{equation*}\label{hatLl}
\widehat {\mathcal L}_\ell=-\frac{d^2}{dr^2} +\frac{(\ell + \frac{d-1}{2}) (\ell+\frac{d-1}{2} -1)}{r^2} + v(r)\,.
\end{equation*}

\begin{proposition}
Let $H_V = -\Delta + V$, where $V(x)= v(r)$ satisfies either  Case A    or Case B   .\\
Any eigenvalue of $-\Delta + V$ is of the form
	\begin{equation}
	\lambda =\lambda_{n,\ell}\,,
	\end{equation}
where $\lambda_{n,\ell}$ is the $n$-th eigenvalue of $\mathcal L_\ell$.\\
A corresponding basis of eigenfunctions has the form
	\begin{equation}\label{unlm}
	u_{n,\ell,m}(r,\omega) = f_{n,\ell} (r) Y_{\ell, m} \,  (\omega)\,,
	\end{equation}
where  $Y_{\ell, m}(\omega)$ denotes an orthonormal basis of (hyper)spherical harmonics. 
\end{proposition}
We recall that these functions form a basis of $L^2(\mathbb{R}^d)$ in Case A     (see \cite{ReSi}) or a basis of the  negative eigenspace  in Case B. 
 \subsection{Courant's nodal  theorem and nodal behavior of eigenfunctions.}\label{ss3.4}
 For the analysis of potentials with singularities it is worth to ask under which condition one can prove Courant's theorem or describe the local nodal structure of an eigenfunction. Under our assumptions the only point
is the control at the origin. Outside the origin, \eqref{ass1} implies that the potential is $C^1$ and the structure of the nodal set is well known.  

Looking at the proof of Courant's theorem, the only thing we need is the unique continuation theorem, i.e. we need to show that if an eigenfunction $u_\lambda $  is identically $0$  in a non empty  open set $\omega$ then it is zero in $\mathbb R^d$. The argument clearly works if there is only a singularity at $0$, because  $\omega \setminus \{0\}$ is an open set where $u_\lambda $ vanishes identically. See \cite{HHOHHO2}  for more properties. We will show in the next subsection that no nodal domain is contained in a sufficiently small ball around the origin. Hence the counting of nodal domains 
  can start outside this ball.

  \subsection{No nodal domains in a small ball}\label{ssfn}
 In this subsection, we show that under our assumptions the nodal domain cannot be contained    in a  small neighborhood of the origin. We will start with Case B which is easier.

 \subsubsection{Case B}
 We have  the following  statement:
\begin{proposition}\label{consHardgen}~\\
If $d\geq 2$ and in Case B there exists $r_d >0$ such that,   if $\lambda < 0$ is an eigenvalue and $u_\lambda$ is a corresponding eigenfunction, there is no nodal domain $\omega$ of  $u_\lambda$  contained in $B(0,r_d)$\,.
  \end{proposition}
  {\bf Proof}\\
  We first  deduce from Assumption \eqref{ass3}, because $s<2$, that: 
  \begin{itemize}
  \item 
  For $d\geq 3$  there exists $r_d$ such that \begin{equation}\label{rd}
  v(r) + \frac{(d-2)^2}{4 r^2} >0\,,
  \end{equation}
   for $r\in (0,r_d)$.  
   \item  For $d=2$, there exists $r_2>0$ such that \begin{equation}\label{r2}
   v(r)  +   \frac{1}{4 r^2\ln^2(r/2 r_2)}  >0\,,
   \end{equation} 
  for $r \in (0,r_2)\,$. 
  \end{itemize}
  We  now use the identity
$$
 \int_\omega |\nabla u_\lambda (x)|^2\, dx  + \int_\omega V(x)  |u_\lambda(x)|^2\, dx  = \lambda \int _\omega |u_\lambda (x)|^2 dx \,,
 $$
 and get because $\lambda < 0$
 $$
 \int_\omega |\nabla u_\lambda (x) |^2\, dx  + \int_\omega V(x) |u_\lambda (x)|^2 \, dx  < 0 \,.
 $$
 When $d\geq 3\,$, we use Hardy's inequality \eqref{Hardyineq}  ($u_\lambda$ is extended in $\mathbb R^d$ by $0$ outside $\omega$ and this extension is in $H^1(\mathbb R^d)$) and get that
 $$
  \int_\omega  \left(V(x)+ \frac {(d-2)^2} {4 r^2}\right)  |u_\lambda (x) |^2 dx  < 0\,.
  $$
  This contradicts \eqref{rd}  if $\omega \subset B(0,r_d)\,$.\\
  When $d = 2\,$, we use the modified Hardy inequality \eqref{logHardy} with $\check R =2 r_2\,$, ($u_\lambda$ is extended by $0$ outside $\omega$ in $D(0,\check R)$ and this extension is in $H_0^1( D(0,\check R)$)) and get a contradiction with \eqref{r2} .

\finbox 

\subsubsection{Case A}
In Case A, with singularities, there is  some difficulty because we consider $\lambda$ large. When $d \geq 3\,$, the previous proposition will be true in a ball whose radius is $r_d(\lambda) \approx \lambda^{-\frac 1 2}\,$. For $d=2\,$, 
$r_2(\lambda)$ could be taken as $r_2(\lambda) \approx \lambda^{-\frac 1 2 -\epsilon}$ for some $\epsilon >0\,$. More precisely, we have
\begin{proposition}
	\label{domainsmallball} Under Assumption  \eqref{ass3}  and if $d\geq 3\,$, 
there exists a constant $c_d>0$ which depends on $V$ and $d$ only and $\lambda_0 >0$  such that, for $\lambda \geq \lambda_0$ and if $u_\lambda$ denotes an eigenfunction of $H_V\,$,  there are no nodal domains of $u_\lambda$  contained in $B(0,c_d\lambda^{-1/2})\,$.\\
If $d =2\,$, for any $\epsilon >0\,$, there exists  $\lambda_\epsilon >0$ and $c_V$ that depends only on $V$  such that, for $\lambda \geq \lambda_\epsilon$ and if $u_\lambda$ denotes an eigenfunction of $H_V$,  there are no nodal domains of $u_\lambda$  contained in $B(0,c_V \lambda^{-1/2 -\epsilon})\,$.
\end{proposition}
\textbf{Proof}\\
	By \eqref{ass3}, there exists $C > 0$ and $r_0 >0\,$, such that $V > -Cr^{-s}$ for \break  $0<r\leq r_0$. As $s\in (0,2)$, there exists $\lambda_0$ such that for $\lambda \geq \lambda_0\,$, there exists  $r_d(\lambda)  \sim\frac{d-2}{2}  \lambda^{-\frac 12} $ such that $$\frac{(d-2)^2}{4} r^{-2} - C r ^{-s} >  \lambda\,,\, \forall r\in (0,r_d(\lambda))\,.
	$$
This implies
$$\frac{(d-2)^2}{4} r^{-2} + v(r) -\lambda  >  0\,,\, \forall r\in (0,r_d(\lambda))\,\mbox{ and } \lambda \geq \lambda_0\,.
$$
	
The proof is achieved by taking $0< c < \frac{d-2}{2}$  in the statement of the proposiiton and  using the Hardy inequality as in the second part of the proof of Proposition \ref{consHardgen}.	\\
For the case $d=2\,$, a not optimal $ r_2(\lambda) =  \lambda^{-\frac 12-\epsilon}$ (for some $\epsilon>0$) together with the modified Hardy inequality do the job for $\lambda_\epsilon$ large enough.\\

\subsection{Upper bound for the degree of the  polynomials associated with the  spherical harmonics }
In this subsection, we prove the existence of a rather optimal upper bound on the degree $\ell$ of the polynomials associated with the  spherical harmonics $Y_{\ell m}$ appearing in the decomposition of an eigenfunction $u_\lambda$.

\begin{proposition}\label{propdegree}~\\
In  Cases A or  B,  if $\lambda$ is an eigenvalue of $H_V$  such that $\lambda < \liminf_{r\ar +\infty}  v$  then  there exists $p_\lambda$  such that, for any  associated eigenfunction $u_\lambda$ and for any $\tau$ satisfying $\inf v < \tau \leq \lambda$, we can find 
	a polynomial $\mathcal P_{\tau,\lambda}$ of $d$ variables of degree at most  $p_\lambda$  such that
	on $V^{(-1)}(\tau)$ in $\mathbb R^d$ the restriction of $u_\lambda$ is equal to the restriction
	of $\mathcal P_{\tau,\lambda}$.\\
	Moreover,  $p_\lambda$ satisfies
	\begin{equation}\label{upbp}
	p_\lambda \leq  \max \{ \ell \, | \,  \ell \geq 1 \mbox{ and } \inf_r \left( v(r) + \frac {\ell (\ell +d-2)}{r^2} \right) < \lambda\ \} \,.
	\end{equation}
\end{proposition}
{\bf Proof.}\\
For given $\lambda$, $u_\lambda$ has the form (see \eqref{unlm})
$$
u_\lambda = \sum_{\lambda = \lambda_{n,\ell}} c_{n,\ell,m} \, u_{n,\ell,m} \,.
$$
If we restrict $u_\lambda$ to the hypersphere of radius $ r_\tau=v^{(-1)} (\tau)$, we get
$$
u_{\lambda}(r_\tau, \omega) =   \sum_{\lambda = \lambda_{n,\ell}} d_{\lambda,\tau,\ell,m}
Y_{\ell,m} (\omega)\,.
$$
Considering the property of $Y_{\ell,m}$, we can choose for the proof of the proposition
\begin{equation}
\mathcal P_{\tau,\lambda}= \sum d_{\lambda,\tau,\ell,m} P_{\ell,m}\,,
\end{equation}
where $P_{\ell,m}$ is the homogeneous  harmonic polynomial of  degree $\ell$ such that  
$$
(P_{\ell,m})_{/\{r=1\}} = Y_{\ell,m}\,.
$$
It remains, in order  to prove \eqref{upbp}, to determine the highest $\ell \geq 1$ such that
$\lambda_{n,\ell} = \lambda$.\\
By the minimax principle, we have 
$$
\lambda=\lambda_{n,\ell}  \geq  m_\ell: =\inf \left( v(r) + \frac { \ell (\ell +d-2)}{r^2}\right) \,.
$$
We have indeed $$\lambda_{n,\ell} = \int_0^{+\infty} |f'_{n,\ell}(r) |^2 r^{d-1} dr + \int_0^{+\infty} \left(v(r) + \frac{\ell (\ell + d-2)}{r^2} \right)| f_{n,\ell}(r)|^2 r^{d-1} dr\,.
$$
\finbox\\~\\
The behavior of $m_\ell$  should be analyzed but note that our assumptions imply that $m_\ell > -\infty\,$. Furthermore $\ell \mapsto m_\ell$ is strictly increasing, so  we can set  $p_\lambda:=[\check p_\lambda] +1 \,$,
where $\check p_\lambda$ is the solution of $\lambda = m_{\check p_\lambda}$ and $[x]$ means the integer part of $x$. \\

{\bf Application: Determination of an upper bound of  $p_\lambda$.}
  
{\bf Case A} \\
We can assume that $\ell \geq 1$. This simply implies later a choice of $p_\lambda \geq 1\,$ 
  If we consider $v(r) =  c \,r^m$ as a model case for $m >1$ and $c>0\,$, the infimum is obtained when
  $$
  c \, m \, r^{m-1} - 2  \frac {\ell (\ell +d-2)}{r^3 } =0\,,
  $$
  i.e. 
  $$
  r = \left(\frac {2 \ell (\ell +d-2)}{c \,m}\right)^{\frac {1}{m+2}} \,.
  $$
  So we get
   $$
   \begin{array}{ll}
  \inf \left( r^m + \frac {\ell (\ell +d-2)}{r^2}\right)  &=\left(\frac {2 \ell (\ell + d-2)}{c \, m}\right)^{\frac {m}{m+2}}
   +  \ell (\ell +d-2) \left(\frac {2 \ell (\ell +d-2)}{c \,m}\right)^{-\frac {2}{m+2}} \\& = \left(\frac{2}{c\,m} \right)^{\frac {m}{m+2}}\, \frac{c\,m+2}{2}\,  (\ell (\ell +d-2))^{\frac {m}{m+2}}\,.
   \end{array}
  $$
  This gives us 
  \begin{equation}
  \check p_\lambda \sim a_m \lambda^{\frac{m +2}{2m} }\,,
   \end{equation}
   with
   $$
   a_m =  (\frac 2 {c\, m})^{-\frac 12}\, (\frac{c \, m+2}{2})^{- \frac{2m}{m+2}}\,.
   $$
  For $m=2\,$, we recover what we got for the harmonic oscillator by direct computation.\\
  To treat the general case, we use the lower bound:
  \begin{equation}\label{ineqcaseA}
  v(r) \geq c \, r^m - \frac{C}{r^s}\,.
  \end{equation}
  We have to estimate
  $$
  \inf \left(  c\, r^m - \frac{C}{r^s}  + \frac { \ell (\ell +d-2)}{r^2}\right) \,.
  $$
  We  observe that:
  $$
   \begin{array}{l}
 \inf \left(  c \, r^m - \frac{C}{r^s}  + \frac { \ell (\ell +d-2)}{r^2}\right)  \\
 \quad \geq   \inf \left(  c \, r^m + \frac { \ell (\ell +d-2)}{2 r^2}\right)  +  \inf \left(   - \frac{C}{r^s}  + \frac { \ell (\ell +d-2)}{2 r^2}\right)  \,.
  \end{array}
$$
  But, there exists (see below the computation in  \eqref{fbelow} with $m=-s$) a constant  $C_0 >0$, such that, for all $\ell \geq 1$, 
  $$
   \inf \left(   - \frac{C}{r^s}  + \frac { \ell (\ell +d-2)}{2 r^2}\right)  \geq - \,C_0\,,
   $$
   and we can use the lower bound of the model case above to get:
   $$
      \inf_r \left( v(r) + \frac {\ell (\ell +d-2)}{r^2} \right) \geq \frac 12\left( \left(\frac{2}{c\, m} \right)^{\frac {m}{m+2}}\, \frac{c\, m+2}{2}\,  (\ell (\ell +d-2))^{\frac {m}{m+2}}\right)   -\, C_0\,.
      $$
   Hence we obtain like for the model case:
   \begin{corollary}\label{corA}
   In Case A, as $\lambda \ar +\infty$\,,
   \begin{equation}\label{asymptplambdaA}
   p_\lambda \approx \lambda^{\frac{m+2}{2 m}} \,.
   \end{equation}
   \end{corollary}
  
    {\bf Case B}\\
  Let us now compute an example corresponding to Case B. 
   If we take $v(r) = - r^m$ for $m \in (-2,0)$, the infimum is obtained when
  $$
  m r^{m-1} + 2  \frac {\ell (\ell +d-2)}{r^3 } =0\,.
  $$
  i.e. 
  \begin{equation}\label{calr}
  r = \left(\frac {2 \ell (\ell +d-2)}{-m}\right)^{\frac {1}{m+2}} \,.
  \end{equation}
  So we get
   \begin{equation}\label{fbelow}
   \begin{array}{ll}
  \inf \left( - r^m + \frac {\ell (\ell +d-2)}{r^2}\right)  &= - \left(\frac {2 \ell (\ell + d-2)}{-m}\right)^{\frac {m}{m+2}}
   +  \ell (\ell +d-2) \left(\frac {2 \ell (\ell +d-2)}{-m}\right)^{-\frac {2}{m+2}} \\& =  - \left(\frac {2}{- m} \right)^{\frac {m}{m+2}}\, \frac{m+2}{2}\,  (\ell (\ell +d-2))^{\frac {m}{m+2}}\,.
   \end{array}
  \end{equation}
  This gives us for $\lambda \ar 0$\,, $\lambda >0$\,,
  \begin{equation}
  \check p_\lambda \sim a_m  (-\lambda)^{\frac{m +2}{2m} }\,,
   \end{equation}
   with
   \begin{equation}
   a_m =  (\frac 2 {-m})^{-\frac 12}\, (\frac{m+2}{2})^{- \frac{2m}{m+2}}\,.
   \end{equation}
   In the Coulomb case $m=-1$ and $d=3$, we get
  \begin{equation}
  \check p_\lambda \sim a_{-1}   (-\lambda)^{-\frac 12}\,,
   \end{equation} 
   to compare with the direct computation which can be done for the Coulomb case.\\
   
   In the general case, we can use
   $$
   v(r) \geq - \,C \, r^{-s}  \,,\, \forall r \in (0,R)
   $$
   and
    $$
    v(r) \geq - \, c \, r^m\,,\, \forall r \in (R,+\infty)\,.
    $$
    
 We will use twice the analysis of the model, the first time with $m$ replaced by $-s$.\\ 
 Let us first consider
 $$
 \inf_{r\in (0,R)}  \left( v(r) + \frac {\ell (\ell +d-2)}{r^2} \right) \geq  \inf_{r\in (0,R)}  \left( - C\, r^{-s}+ \frac {\ell (\ell +d-2)}{r^2} \right)\,.
 $$
 We observe (see \eqref{calr} with $m=-s$)  that for $\ell$ large enough the map \break   $r \mapsto  - C\, r^{-s}+ \frac {\ell (\ell +d-2)}{r^2}$ is decreasing on $(0,R)$. Hence
 $$
  \inf_{r\in (0,R)}  \left( - C\,  r^{-s}+ \frac {\ell (\ell +d-2)}{r^2} \right)=    \left( - C\, R^{-s}+ \frac {\ell (\ell +d-2)}{R^2} \right) \geq - C_s\,.
  $$
  For the second case, we can use 
  $$
  \inf_{r\in (R,+\infty)}  \left(\, - \, C\, r^{m}+ \frac {\ell (\ell +d-2)}{r^2} \right)\geq  \inf_{r\in (0,+\infty)}  \left( - C\, r^{m}+ \frac {\ell (\ell +d-2)}{r^2} \right)\,,
  $$
  and what we obtained for the homogeneous model.\\
  
  \begin{corollary}\label{corB}
   In Case B, 
    \begin{equation}\label{asymptplambdaB}
   p_\lambda \approx (- \lambda) ^{\frac{m+2}{2m}} \,.
   \end{equation}
   \end{corollary}

\subsection{Nodal domains on hyperspheres}\label{ss3.7}

Since the considered potentials $V$  are radial, the energy hypersurfaces  $\left\lbrace  V= \alpha \lambda \right\rbrace $ are hyperspheres centered at $0$. Also, the restriction of any eigenfunction $u_\lambda $ of $H_V$ to a hypersphere equals the restriction of some  harmonic polynomial.
We can use the following result proven in \cite{Cha1}, which is based on  \cite{Mi}:

 \begin{proposition}
 	\label{lemmapoly1}
 	Let $P$  be a polynomial of degree $k$ with $d$ variables. Then its restriction to the hypersphere $\mathbb S^{d-1}$ admits at most $2^{2d-1} k ^{d-1}$ nodal domains.
 \end{proposition}

We will combine this with the previous estimates obtained in Corollaries~\ref{corA} and \ref{corB}  to obtain an upper bound on the number of nodal domains on any hypersphere.

 \section{Weyl's formula}\label{s4}
\subsection{Preliminaries}
For Schr\"odinger operators, Weyl's formula takes (under of course suitable assumptions to be discussed below)  the form
\begin{equation}
N (\lambda)  \sim (2\pi)^{-d} \int_{\xi^2 + V(x) \leq \lambda} dx d\xi\,.
\end{equation}
After integration in the $\xi$ variable, we get
\begin{equation}\label{Wformula}
N (\lambda)  \sim W(\lambda) \,,
\end{equation}
where 
\begin{equation}
W(\lambda):= w_d \int (\lambda - V)^{\frac d2}_+ dx\,,
\end{equation}
with $w_d$ defined in \eqref{defw}.\\

This formula makes sense  in case A (as $\lambda \ar +\infty$) and  in case B (as $\lambda \ar 0$ with $\lambda < 0$). Let us just compute the right hand side for the two toy models: the harmonic oscillator and the Schr\"odinger operator with Coulomb potential.
For the harmonic oscillator, we simply get
\begin{equation}\label{asharm}
W(\lambda)= w_d \int (\lambda - r^2)^{\frac d2}_+ dx =  h_d\, w_d \, \lambda^d \,,
\end{equation}
with
$$
h_d:=  \int (1 - r^2)^{\frac d2}_+ dx >0\,.
$$
More generally, if $v(r)= r^m$ for $m>0$, we obtain, as $\lambda \ar +\infty$,
\begin{equation}
W (\lambda) = w_d\, h_{d,m}\, \lambda^{d (\frac12 +\frac 1m)}\,.
\end{equation}
In the Coulomb case, we get, with $\lambda < 0$
\begin{equation}
W(\lambda) =  w_d   \int \left(\lambda + \frac 1 r\right)^{\frac d2}_+ dx = e_d\, w_d  \,  (-\lambda)^{- \frac d 2} \,,
\end{equation}
with
 \begin{equation}
 e_d:=\int \left(\frac{1}{|x|} -1\right)^{\frac d 2}_+ dx = |\mathbb S^{d-1}| \int_0^1  \left(1-r\right)^{\frac d 2} r^{\frac d2 -1} dr < +\infty\,.
 \end{equation}
 More generally, if $v(r) =- r^m$, for $m\in (-2,0)$, we obtain as $\lambda \ar 0$ ($\lambda < 0$),
 \begin{equation}
W(\lambda)  = w_d h_{d,m}  |\lambda|^{d (\frac12 +\frac 1m)}\,.
 \end{equation}
Observing that $N(\lambda_n) =n-1$ if $\lambda_{n-1} < \lambda_n$, and assuming that the Weyl formula is proven (see below for the proof), we get conversly 
\begin{equation}
\lambda_n \sim \check w_d\, n^{\frac 1d}\quad  \mbox{ with } 1=h_d\, w_d \,(\check w_d)^d\,,
\end{equation}
in the case of the harmonic oscillator and
\begin{equation}\label{WeCo}
\lambda_n \sim - \check v_3 \, n^{-\frac 23}\,,\mbox{ with }\quad  1 = e_3 \,w_3 \, (\check v_3)^{-\frac 32}\,,
\end{equation}
in the case of the Coulomb case.\\~\\

More generally we have the proposition:
\begin{proposition}
In  Cases A or  B
\begin{equation}\label{Asw}
W(\lambda) \approx |\lambda|^{d (\frac12 +\frac 1m)}\,,
\end{equation}
where the asymptotics is as $\lambda \ar +\infty$ in Case A and as $\lambda \ar 0$ ($\lambda < 0$) in Case B.
\end{proposition}
{\bf Proof}\\
Outside a ball we can use for estimating  the integral defining $W(\lambda)$ the comparison of $v(r)$ with $r^m$ and then use the previous computations for the models. The control of the integral in a ball will be done in Subsection \ref{ss4.3}.
\finbox

\subsection{Weyl's formula under weak assumptions}
There is  vast literature on this subject: Reed-Simon \cite{ReSi} and references therein (for the historics), D.~Robert \cite{Rob}, H. Tamura \cite{Ta0}, Tulovski-Shubin \cite{TS},  R. Beals \cite{Be}, L. H\"ormander \cite{Ho,Ho1},  A. Mohamed \cite{Mo}. In  the recent contributions  the goal is to control the remainder but this is  not important in the applications considered here.
Here, we prefer to work under weaker asssumptions and  can use Theorem XIII.81 in Reed-Simon (Vol. 4) \cite{ReSi} for the case $V\ar +\infty$ with a condition $d\geq 2$ and $m >1$,  and  for the case $V\ar 0$ as $|x| \ar + \infty$,  Theorem XIII.82. The treatment of the singularity is also explained (without detail) (see the discussion p. 277, lines -7 to -1, sending to  Problem 132 therein). The idea there is to first prove a statement with $V$ continuous  and then to show that the addition of a potential $W$ with compact support or in $L^\frac d2$  ($d\geq 3$) does not change the Weyl asymptotics. \\

Theorem XIII.81 in \cite{ReSi} reads:
\begin{theorem}\label{WRS} Let $V$ be a measurable function on $\mathbb R^d$ ($d\geq 2$) obeying
\begin{equation} 
c_1 \, (r^\beta -1) \leq V(x) \leq c_2 \,(r^\beta +1)\,,
\end{equation}
and
\begin{equation}\label{lip}
| V(x)-V(y)| \leq c_3 \, [ \max \{ |x|,|y|\}]^{\beta -1} |x-y|\,,
\end{equation}
for some $\beta >1$ and suitable constants $c_1,\,c_2,\,c_3 >0\,$.\\
Then 
$$
\lim_{\lambda \ar +\infty} N (\lambda) / W(\lambda) =1\,.
$$
\end{theorem}

\begin{remark}
The theorem is still true if we consider the Dirichlet problem for $H_V$ in $\mathbb R^d \setminus B$, where $B$ is a ball centered at $0$. In Case A, the assumptions of the theorem are satisfied in $\mathbb R^d \setminus B$.
This follows of our assumption \eqref{ass5} on $v'$.
\end{remark}

For Case B,  Theorem XIII.82 in \cite{ReSi} reads:
\begin{theorem} Let $V$ be a measurable function on $\mathbb R^d$ ($d\geq 2$) obeying
\begin{equation} - c_1\, (r +1)^{-\beta} \leq V(x) \leq  - c_2\, (r +1)^{-\beta}  \,,
\end{equation}
and
\begin{equation}
| V(x)-V(y)| \leq c_3\, [ 1+ \min  \{ |x|,|y|\}]^{-\beta -1} |x-y|\,,
\end{equation}
for some $\beta <2$ and suitable constants $c_1,\,c_2,\,c_3 >0\,$.\\
Then 
$$
\lim_{\lambda \ar +\infty} N (\lambda) / W(\lambda) =1\,.
$$
\end{theorem}
\begin{remark}
The theorem is still true if we consider the Dirichlet problem for $H_V$ in $\mathbb R^d \setminus B$, where $B$ is a ball centered at $0$. In Case B, the assumptions of the theorem are satisfied in $\mathbb R^d \setminus B$. This is a consequence of our assumption on $v'$ in \eqref{ass7}.
\end{remark}

  \subsection{Treatment of the singularity}\label{ss4.3}
 To cover the question of the treatment of the singularity at the origin we could think of  using (Problem 132 in \cite{ReSi}) to treat the singularity as a perturbation. Due to the use of the Cwickel-Lieb-Rozenblum inequality \cite{Cw,Li,Roz} in the argument, this approach works only under the condition  $d\geq 3\,$.  If we remember that we only need a lower bound for $N (\lambda) $ one can proceed for $d\geq 2$ in the following way.  We can introduce a small ball $B=B(0,\epsilon) $ around the singularity and look at the Dirichlet problem in $\mathbb R^d \setminus B\,$. 
  We denote by $N_B(\lambda)$ the corresponding counting function.   
  Because the eigenvalues are greater than the initial problem by monotonicity of the domain, the estimate of the  $N( \lambda)$ of the new problem will give the lower bound:
  $$
  N_B(\lambda) \leq N (\lambda) \,.
  $$
  The theorem in Reed-Simon \cite{ReSi} can be applied in $\mathbb R^d \setminus B$ (proof unchanged) and we get by Weyl's formula
  $$
  N_B(\lambda) \sim W_B(\lambda)\,,
  $$
  with
  $$
  W_B(\lambda)= w_d\, \int_{r\geq \epsilon} (\lambda -V)_+^{\frac d 2}\, dx\,.
  $$
  It remains to compare $W_B(\lambda)$ and $W (\lambda)$ in our two cases.\\

{\bf  Case B}\\
 Here $\lambda < 0\,$. It is enough  to show that, for some $\epsilon >0$,
 $$ 
 \int_{B(0,\epsilon)} (\lambda- V)_+^{\frac d2} \, dx  < +\infty\,.
$$
  We have $$\int_{B(0,\epsilon)} (\lambda- V)_+ ^{\frac d2} \,dx\leq C\, \int_0^\epsilon  r^{d-1-s \frac d2} \, dr < +\infty\,,$$ the finiteness resulting from  the assumption $s <2\,$.  \\

{\bf Case A}.\\
 Here $\lambda \geq \lambda_0 >0$.  We will show that $\int_{B(0,\epsilon)} (\lambda- V)_+ ^{\frac d2} \,dx$ is relatively small in comparison with $N(\lambda) $. In   Case A, we have seen that 
\begin{equation}\label{asanharm}
W(\lambda)  \approx  \int (\lambda - r^m)^{\frac d2}_+\, r^{d-1} dr  \approx \lambda^{\frac d2 + \frac dm} \,,
\end{equation}
We have 
$$
\int_{B(0,\epsilon)} (\lambda- V)_+ ^{\frac d2} \,dx\leq C \int_0^\epsilon (\lambda + r^{-s})^{\frac d2} \, r^{d-1} dr \leq \hat C \, (\hat C\, +\,\lambda^\frac d2) \,.
$$
This gives, as $\lambda \ar +\infty\,$,
$$
\int_{B(0,\epsilon)} (\lambda- V)_+ ^{\frac d2} \,dx/ W (\lambda)  = \mathcal O (\lambda^{-\frac dm})\,.
$$

Hence in the two cases, we have shown that $W_A(\lambda) \sim W (\lambda)\,$. In conclusion, we have obtained the following theorem.
\begin{theorem}
	\label{theoremeweylfinal}~\\
In  Cases A or  B\,, if $d\geq 2\,$, we have
\begin{equation}\label{AsWN}
N (\lambda)  \geq W (\lambda) (1 + o (1))\,,
\end{equation}
where the remainder  $o(1)$ is as $\lambda \ar +\infty$ in Case A and as $\lambda \ar 0$ ($\lambda <0$) in Case B\,.
\end{theorem}

\section{Counting nodal domains}\label{s5}
\subsection{Preliminaries}
We construct a radial partition of $\{V < \lambda\}$ of cardinality $\nu (\lambda) $ with $\nu (\lambda) $ to be defined later.
When $v$ is increasing on $(0,+\infty)$, $r_\lambda:=v^{-1} (\lambda)$ is well defined  in $(\inf v,+\infty)$ in Case A, and for any $\lambda$ in   $(-\infty,0)$ in Case B.  But we will only be  interested  $\lambda \ar +\infty$ in Case A, and in $\lambda \ar 0$ ($\lambda <0$) in Case~B.\\ Under the weaker Assumption \eqref{assmonotone}, we can define in the two cases  $r_\lambda$ by 
$$
r_\lambda := \sup_{v(r)=\lambda} r\,,
$$
and obtain the asymptotics 
$$
r _\lambda \approx  |\lambda|^{\frac 1m}\,.
$$
Note that in both cases $r_\lambda$ tends to $+\infty\,$.
\subsection{Analysis of Case A}
\subsubsection{A first partition of the classical region.}

We recall that in this case, we assume
$
\lambda \geq \lambda_0 >0\,.
$

We introduce a partition of the classical region by introducing $\nu (\lambda) $  annuli, for $i=1,\dots,\nu (\lambda) $,
$$
D_i(\lambda):= \left\{ x \in \mathbb R^d \,|\, \left( \frac{i-1}{\nu (\lambda) }\right)^{1/d}r_\lambda  < r < \left(\frac{i}{\nu (\lambda) }\right)^{1/d}r_\lambda \right\}\,.
$$
We note that each annulus has the same volume:
\begin{equation}\label{D2}
|D_i(\lambda)| = \omega_d  \frac{1}{\nu (\lambda) } {r_\lambda }^d \approx \frac{\lambda^{\frac dm}}{\nu (\lambda) }\,.
\end{equation}
The cardinality $\nu (\lambda) $ satisfies a priori the condition
\begin{equation}\label{condn1}
\lim_{\lambda \ar +\infty} \nu (\lambda)  = +\infty\,,
\end{equation}
but  the condition
\begin{equation}\label{condn2}
\lim_{\lambda \ar +\infty} \nu (\lambda) ^{-1/d} \, r_\lambda = +\infty\,,
\end{equation}
 will appear along the proof.\\

The determination of $\lambda_0$ "large enough"  will be given during the proof.
If $u_\lambda$ denotes some eigenfunction, we denote by 
$\mathcal D (u_\lambda)$ the set of the nodal domains of $u_\lambda$. We now introduce in $\mathcal D (u_\lambda)$ the following subsets.

\begin{definition}
	\label{ai}
	$$A_i(u_\lambda)= \left\{ \Omega \in \mathcal D (u_\lambda) \,|\, \Omega \subset D_i(\lambda) \right\}\, .$$ Here, $i$ can take the values $1,2, \ldots, \nu (\lambda) $.
\end{definition}

\begin{definition}
	\label{bj}
	$$B_j(u_\lambda) = \left\{ \Omega \in \mathcal D(u_\lambda)\, \left| \, \, \Omega \cap \left\lbrace r = \left({\frac{j}{\nu (\lambda) }}\right)^{1/d}r_\lambda  \right\rbrace \neq \emptyset   \right. \right\} \, .$$ Here, $j$ can take the values $1,2, \ldots, \nu (\lambda) \,$.
\end{definition}

 From Subsection \ref{ss3.2}, we know that every nodal domain is contained in at least one of these sets.\\

\begin{remark}
 This partition will be refined by introducing, in the case of a singularity at the origin,  a  further partition of $A_1(u_\lambda)$.
 \end{remark}

\subsubsection{Counting the nodal domains contained in one annulus of the partition}
We first count in each of the annuli $D_i(\lambda)$ for $i\geq 2\,$.  Except if there are no singularity, the treatment of $A_1(u_\lambda)$ will be done separately. Hence we first prove the
\begin{proposition}\label{Prop5.4}~\\
In Case A, if $\nu (\lambda) $ satisfies \eqref{condn1} and \eqref{condn2}, 
we have the following inequality,  as $\lambda \ar +\infty$ , 
\begin{equation}\label{ineqcomp1} 
\sum_{i=2}^{\nu (\lambda) } \# A_i(u_\lambda) \leq \gamma(d) \, W(\lambda) \, \left(1 +  \mathcal O \left(\frac{1}{\nu (\lambda) }\right)\right)\,.
	\end{equation}
\end{proposition}
{\bf Proof}\\
If  $\Omega$ is a bounded nodal domain of $u_\lambda \,$, the  Faber-Krahn inequality (see Theorem~\ref{ThFK}) gives:
\begin{equation}
\label{ineqaire3}
\frac{\int_\Omega {{|\nabla u_\lambda (x) |}^2}\,dx }{\int_\Omega {u_\lambda(x) ^2}\,dx } \geq \left(\frac{1}{|\Omega|}\right)^{\frac{2}{d}}{ \omega_d }^{\frac{2}{d}} \lambda (B_d) \, ,
\end{equation}
where we recall that $|\Omega|$ denotes the volume of $\Omega\,$.\\
We know that for a given bounded nodal domain $\Omega$, we have
\begin{align*}
\lambda \, = \frac{\int_\Omega {{|\nabla u_\lambda |}^2}dx + \int_\Omega {V u_\lambda ^2 }\,dx }{\int_\Omega {u_\lambda ^2}\,dx} \, ,
\end{align*}
which implies
\begin{equation}\label{inepaireomega}
\frac{\int_\Omega {{|\nabla u_\lambda |}^2} \,dx}{\int_\Omega {u_\lambda ^2}\,dx } < \lambda -  \inf_{x\in \Omega} V(x)  \, .
\end{equation}
For all $\Omega \in  A_i(u_\lambda)$, we obtain that 
\begin{align} \label{ineqaire}
\frac{\int_\Omega {{|\nabla u_\lambda |}^2} \, dx}{\int_\Omega {u_\lambda ^2}\, dx } < \lambda - v\left(\left({\frac{i-1}{\nu (\lambda) }}\right)^{1/d}r_\lambda \right) \, .
\end{align}
Because $i \geq 2\,$, we can, under Condition \eqref{condn2}, 
combine \eqref{ineqaire3} and \eqref{ineqaire} and  obtain
\begin{equation} \label{ineqvolume}
|\Omega| \geq  \frac{ \omega_d \,  \lambda (B_d) ^{\frac d2}}{\left(\lambda - v\left(\left(\frac{i-1}{\nu (\lambda) }\right)^{1/d}r_\lambda \right)\right)^{\frac{d}{2}}} \, .
\end{equation}
Observing that $$\sum_{\Omega \in A_i(u_\lambda)}|\Omega| \leq |D_i(\lambda)|\,,$$
 we obtain that
\begin{equation}\label{contint}
 \mathrm{\# }\,A_i(u_\lambda)\, \leq \frac{1}{ \omega_d \lambda(B_d)^{\frac d2}}\,\left(  |D_i(\lambda)|   \left(\lambda - v\left(\left(\frac{i-1}{\nu (\lambda) }\right)^{1/d}r_\lambda \right)\right)^{\frac{d}{2}}\right) \,.
\end{equation}
Summing up for $i=2,\dots,\nu (\lambda) $, we get
\begin{align}
\sum_{i=2}^{\nu (\lambda) } \mathrm{\# }\,A_i(u_\lambda)\, \leq \frac{1}{ \omega_d \lambda(B_d)^{\frac d2}}\,\left( \sum_{i=2}^{\nu (\lambda) }   |D_i(\lambda)|   \left(\lambda - v\left(\left(\frac{i-1}{\nu (\lambda) }\right)^{1/d}r_\lambda \right)\right)^{\frac{d}{2}}\right) \,.
\end{align}
We recognize on the right hand side a  Riemann sum for the function \break  $x \mapsto (\lambda-V(x))^{\frac d2}$ in $D(0, r_\lambda) \setminus (D_{\nu (\lambda) }(\lambda) \cup D_1(\lambda))$.\\
More precisely, we can write, using the monotonicity of $v\,$,
$$
\begin{array}{l}
 \sum_{i=2}^{\nu (\lambda) }   |D_i(\lambda)|   \left(\lambda - v\left(\left(\frac{i-1}{\nu (\lambda) }\right)^{1/d}r_\lambda \right)\right)^{\frac{d}{2}}\\
 \hspace{3cm}\leq |D_2(\lambda)|\,   \left(\lambda - v\left(\left(\frac{1}{\nu (\lambda) }\right)^{1/d}r_\lambda \right)\right)^{\frac{d}{2}}  + \int (\lambda- V(x))_+ ^\frac d2\, dx\\
 \hspace{3cm} \leq |D_2(\lambda)| \, \lambda^{\frac d2} + \int (\lambda- V(x))_+ ^\frac d2 \, dx \,.
 \end{array}
$$

Using the asymptotic behavior of $v$ at $+\infty$ given in  \eqref{ass6}, the computation of  $\gamma(d)$ in \eqref{gamma2}, the asymptotic behavior  of $W(\lambda)$ given in \eqref{Asw},  and \eqref{D2}, we achieve the proof of
 Proposition \ref{Prop5.4}.
\subsubsection{Counting the nodal sets meeting the boundary of the annuli}
Let us now turn to the study of the sets $B_i(u_\lambda)$. We have shown in Corollary~\ref{corA} that $p_\lambda \approx \lambda^{\frac{m+2}{2m}}$.  Using Proposition \ref{lemmapoly1}, we obtain 
 that the number of nodal domains in a given $B_i(u_\lambda)$ satisfies
\begin{equation}\label{zrobd}
\textrm{\#}B_i(u_\lambda) \leq 2^{2d-1}\,  p_\lambda^{d-1} \leq C_d\, \lambda^{\frac{(d-1)(m+2)}{2m}}\,.
\end{equation}

Comparing with \eqref{Asw}, we get that, for some $C>0$ and $\lambda \geq \lambda_0$\,,
\begin{equation*}
\sum_i \textrm{\#}B_i (u_\lambda) \leq  C \,  \nu (\lambda) \, \lambda^{-\frac{m+2}{2m}} \, W(\lambda)  \,.
\end{equation*}
If $\nu (\lambda) $ satisfies in addition, 
\begin{equation}\label{condn3}
\lim_{\lambda \ar +\infty} \nu (\lambda) \,  \lambda^{-\frac{m+2}{2m}} =0\,, 
\end{equation}
we obtain
\begin{equation}\label{Asyb}
\lim\limits_{\lambda \to +\infty}\frac{\sum_i  \, \#B_i(u_\lambda)}{W(\lambda)} = 0 \,.
\end{equation}

\subsubsection{Counting in $D_1(\lambda)$}

We still have to consider when there is a singularity around $0$. We treat first the case $d \geq 3\,$.
 We know from Proposition \ref{domainsmallball}  that we can replace $D_1(\lambda)$ by
 the annulus $\widehat D_{1}(\lambda,C)$ defined by
   $$
  \widehat D_{1}(\lambda,C):=\{x\in \mathbb R^d\,|\,  \frac 1C \lambda^{-\frac 12} < r < r_\lambda/\nu (\lambda)^{\frac{1}{d}} \}$$
   for a sufficiently large $C$. \\
   The number of nodal domains $\mu_{10} (\lambda)$ crossing the hypersphere  $\{ r= \frac 1C \lambda^{-\frac 12} \}$  is controlled by \eqref{zrobd}:
   \begin{equation}\label{mu01}
   \lim_{\lambda \ar +\infty} \frac{\mu_{10} (\lambda)}{W(\lambda)} =0\,.
   \end{equation}
   
 To continue, we consider a partition of $\widehat D_1(\lambda)$ in two annuli:
   $$
 D_{11}(\lambda) := \{x\in \mathbb R^d\,|\,  \frac 1C \lambda^{-\frac 12} < r < C \}\mbox{ and } D_{12}(\lambda):=\{x\in \mathbb R^d\,|\,  C < r < r_\lambda/\nu (\lambda)^\frac 1d  \}\,,
   $$
   where we keep the liberty to choose $C  >0$ larger than the previous one.\\
  Again the   number $\mu_{11}(\lambda)$  of nodal domains crossing the hypersphere  $\{ r= C \}$  is controlled by \eqref{zrobd}:
    \begin{equation}\label{mu11}
   \lim_{\lambda \ar +\infty} \frac{\mu_{11} (\lambda)}{W(\lambda)} =0\,.
   \end{equation}

  {\bf Control in $D_{12}(\lambda)$.\\}
 The treatment of $D_{12}(\lambda)$ can be done for $C$ large enough (in order to have the monotonicity of $v$) like the analysis of the $D_i(\lambda)$ for $i\geq 1$. More precisely, we can replace \eqref{contint} (for $i=1$) by
 \begin{equation}
 \mathrm{\# }\,A_{12}(u_\lambda)\, \leq \frac{1}{ \omega_d \lambda(B_d)^{\frac d2}}\,   |D_{12} (\lambda)|  \,  \left(\lambda - v(C)\right)^{\frac{d}{2}}\,,
\end{equation}
 where  
 $$A_{12}(u_\lambda):= \{ \Omega \in \mathcal D(u_\lambda)\,|\, \Omega \subset D_{12} (\lambda) \} \,.$$
 Hence we get, for some constant $C_d>0\,$,
  \begin{equation}
 \mathrm{\# }\,A_{12}(u_\lambda)\, \leq C_d \,\nu(\lambda)^{-1} \, \lambda^{\frac d2 + \frac dm}\,,
\end{equation}
which implies
\begin{equation}\label{contd12}
\lim_{\lambda \ar + \infty} \left(\mathrm{\# }\,A_{12}(u_\lambda)/ W(\lambda)\right)=0\,.
\end{equation}
 {\bf Control in $D_{11}(\lambda)$.\\}
  Note that in $D_{11}(\lambda)$, we have, for some constant $C_s >0\,$, 
  $$
  V(x) \geq - C_s \, \lambda^{\frac{s}{2}}\,,\, \forall x \in D_{11}(\lambda)\,.
  $$
  Hence, for $\lambda \geq \lambda_0$,  we deduce from \eqref{inepaireomega}, that 
  $$
  \frac{\int_\Omega |\nabla u_\lambda(x)|^2\,dx }{\int_\Omega u_\lambda(x)^2\,dx} \leq 2 \lambda \,.
  $$
  As in the proof of \eqref{ineqvolume} we obtain:
  \begin{equation}
  | \Omega| \geq  \omega_d \,\lambda(B_d)^\frac d2 \, 2^{-\frac d2} \, \lambda^{-\frac d2}\,.
  \end{equation}
  Summing over the $\Omega$'s contained in $D_{11}(\lambda)$ and observing that the volume of $D_{11}(\lambda)$ is bounded, we get the existence of a constant $\hat C_d$ such that
  $$
  \# A_{11}(u_\lambda)\leq \hat C_d \, \lambda^{\frac d2}\,,
  $$
  where 
  $$A_{11}(u_\lambda):= \{ \Omega \in \mathcal D(u_\lambda)\,|\, \Omega \subset D_{11} (\lambda) \} \,.$$
  In particular we get
  \begin{equation}\label{lastasy}
    \lim_{\lambda \ar + \infty} \frac{ \# A_{11}(u_\lambda)
 }{W(\lambda)} =0\,.
  \end{equation}
  The case when $d=2$ does not lead to new difficulties.
  \subsubsection{Conclusion for Case A}
  Summing all the upper bounds and having chosen $\nu (\lambda) $ satisfying \eqref{condn1}, \eqref{condn2}, and \eqref{condn3},  we get, as $\lambda \ar +\infty\,$,
  \begin{equation}
  \mu(u_\lambda) \leq  \gamma (d)  \, W(\lambda) \, (1 + o(1))\,.
\end{equation}
Using the asymptotic upper bound  \eqref{AsWN}, we get, as $\lambda \ar +\infty\,$,
\begin{equation}
\label{eqmu1}
  \mu(u_\lambda) \leq  \gamma (d)  \, N (\lambda)  \, (1 + o(1))\,.
\end{equation}
Observing that $N(\lambda_n) \leq  n -1\,$, we obtain  Theorem \ref{thm1.6}.

\subsection{Case B}
The proof in case B follows the same lines.
We define the sets $D_i(\lambda)$, $A_i(u_\lambda)$ and $B_i(u_\lambda)$ for $i=2, \ldots, \nu(\lambda)$ as in case A, with $\nu(\lambda)$ satisfying conditions equivalent to  (\ref{condn1}), (\ref{condn2}) and (\ref{condn3}), replacing $\lambda \ar \infty$ by $\lambda \ar 0$ (with $\lambda <0$). Hence, we assume
that $\nu (\lambda) $ satisfies the conditions
\begin{equation}\label{condn1B}
\lim_{\lambda \ar 0} \nu (\lambda)  = +\infty\,,
\end{equation}
\begin{equation}\label{condn2B}
\lim_{\lambda \ar 0} \nu (\lambda) ^{-1/d} \, r_\lambda = +\infty\,,
\end{equation}
and
\begin{equation}\label{condn3B}
\lim_{\lambda \ar +0} \nu (\lambda) \,  |\lambda|^{-\frac{m+2}{2m}} =0\,.
\end{equation}

\begin{proposition}~\\
In Case B, if $\nu (\lambda) $ satisfies \eqref{condn1B} and  \eqref{condn2B}, then
we have the following inequality,  as $\lambda \ar 0$  ($\lambda <0\,$), 
\begin{equation}\label{ineqcomp1}  
\sum_{i=2}^{\nu (\lambda) } \# A_i(u_\lambda) \leq \gamma(d) \, W(\lambda) \, (1 +  \mathcal O (\frac{1}{\nu (\lambda) }))\,.
\end{equation}
\end{proposition}
The proof is the same as in case A. For $W(\lambda)$, we can use  the asymptotics \eqref{Asw}.\\

 For the cardinalities of the sets $B_i(u_\lambda)$,  (\ref{Asyb}) holds in case B under condition \eqref{condn3B} and we can use  Corollary \ref{corB} together with \eqref{Asw}.\\

The treatment of the singularity is slightly easier in this case. We use Proposition \ref{consHardgen} to make a partition of $D_1(\lambda)$ in  two annuli:
$$
D_{11}(\lambda,C) := \{x\in \mathbb R^d\,|\,  r_d < r < C \} \,,
$$
and 
$$
 D_{12}(\lambda,C):=\{x\in \mathbb R^d\,|\,  C < r < r_\lambda/\nu (\lambda)^\frac 1d  \}\,.
$$

Again, we choose $C$ such that $v$ is strictly increasing for $|x| > C$. 
Since $\lambda < 0$, there exists $M > 0 $ such that $\# A_{11}(u_\lambda) < M\,$. 
To give an upper bound on $\# A_{12}(u_\lambda)$, we follow the same steps as in case A.\\

The behavior of the number of eigenvalues below or equal to $\lambda$ is given by  (\ref{AsWN}), hence inequality (\ref{eqmu1}) is verified.
Since $N(\lambda_n) \leq n-1\,$, we obtain Theorem \ref{thm1.6} in this case as well.


\begin{thebibliography}{99}

\bibitem{Be} R. Beals.
\newblock A general calculus of pseudo-differential operators.
\newblock Duke Math. J. 42 (1975), 1--42.

\bibitem{BeHe} P. B\'erard, B. Helffer.
\newblock     On the nodal patterns of the 2D isotropic quantum harmonic oscillator.
\newblock arXiv:1506.02374 (2015).


\bibitem {BeMe} P. B\'erard, D. Meyer.
\newblock In\'egalit\'es isop\'erim\'etriques et applications.
\newblock Annales scientifiques de l'\'Ecole Normale Sup\'erieure, S\'er. 4, 15(3) (1982), 513-541.

\bibitem{Bo} J. Bourgain.
\newblock On Pleijel's nodal domain theorem. 
\newblock ArXiv:1308.4422.  Int. Math. Res. Not. (2015), no. 6, 1601--1612.

\bibitem{CaHu} L. De Carli, S.M. Hudson.
\newblock A Faber-Krahn inequality for solutions of
Schr\"odinger's equation.
\newblock 
 Advances in Mathematics 230 (2012), 2416--2427.

\bibitem{Cha} P. Charron.
\newblock Th\'eor\`eme de Pleijel pour l'oscillateur
harmonique quantique.
\newblock M\'emoire de Ma\^\i trise (2015). Universit\'e de Montr\'eal.


\bibitem{Cha1} P. Charron.
\newblock  A Pleijel-type theorem for the quantum harmonic oscillator.
\newblock  arXiv:1512.07880. To appear in Journal of Spectral Theory (2016).

\bibitem {Cou} R. Courant.
\newblock Ein allgemeiner Satz zur Theorie der Eigenfunktionen selbstadjungierter Differentialausdr\"ucke.
\newblock Nachr. Ges. G\"ottingen (1923), 81-84.


\bibitem {CH} R. Courant and D. Hilbert.
\newblock Methods of Mathematical Physics, Vol. 1.
\newblock Wiley New York (1953).

\bibitem{Cow} C. Cowan.
\newblock Optimal Hardy inequalities for general elliptic operators with improvements.
\newblock Comm. Pure Appl. Math. 9 (2010), 109--140.  

\bibitem{Cw} M. Cwikel. 
\newblock Weak type estimates for singular values and the number of 
bound states.
\newblock Annals of Math. 106 (1977), 93--100.

\bibitem {Fab} C. Faber.
\newblock Beweis, dass unter allen homogenen Membrane von gleicher
Fl\"ache und gleicher Spannung die kreisf\"ormige die tiefsten Grundton
gibt.
\newblock Sitzungsber. Bayer. Akad. Wiss.,
Math. Phys. M\"{u}nchen (1923), 169-172.

\bibitem{Hel} B. Helffer.
\newblock  Spectral theory and applications.
\newblock Cambridge University Press (2013).

\bibitem{HKSW} B. Helffer, A. Knauf, H. Siedentop, and R. Weikard.
\newblock On the absence of a first order correction of a Schr\"odinger operator with Coulomb singularity.
\newblock CPDE 17 (3\&4) (1992), 615-639.

\bibitem{HHO} B. Helffer and T. Hoffmann-Ostenhof.
\newblock A review on large k minimal spectral $k$-partitions and Pleijel's Theorem. 
\newblock ArXiv:1509.04501. Spectral theory and partial differential equations, 39--57  (2015),  Contemp. Math. 640.



\bibitem{HPS} B. Helffer and M. Persson-Sundqvist.
\newblock  On nodal domains in Euclidean balls.
\newblock ArXiv:1506.04033v2. To appear in Proc. AMS 2016.


\bibitem{HHOHHO2} M. and T. Hoffmann-Ostenhof.
\newblock Local properties of solutions of Schr\"odinger equations.
\newblock Comm. in Partial Differntial Equations  17 (3\&4) (1992), 491--522.

\bibitem{Ho} L. H\"ormander.
\newblock The Weyl calculus of pseudodifferential operators.
\newblock Comm. in Pure Appl. Math. 32, 359--443 (1979).

\bibitem{Ho1} L. H\"ormander.
\newblock  On the asymptotic distribution of the eigenvalues of pseudodifferential operators in $\mathbb R^n$.
\newblock  Arkiv f\"or matematik 17 (3)  297--313 (1979).


\bibitem {Kr} E. Krahn.
\newblock \"Uber eine von Rayleigh formulierte minimal Eigenschaft
des Kreises.
\newblock Math. Ann. 94 (1925), 97-100.

\bibitem{Ley} J. Leydold.
\newblock Knotenlinien und Knotengebiete von Eigenfunktionen. 
\newblock Diplom Arbeit, Universit\"at Wien (1989).

\bibitem{Li} E. H. Lieb.
\newblock  Bounds on the eigenvalues of the Laplace and Schr�dinger 
operators.
\newblock Bull. Am. Math. Soc. 82 (1976), 751--753.


\bibitem{Mi} J. Milnor.
\newblock  On the Betti Numbers of Real Varieties. 
\newblock Proc. Amer. Math. Soc. 15 (1964), 275-280.

\bibitem{Mo} A. Mohamed.
\newblock  Comportement asymptotique, avec estimation du reste, des valeurs propres d'une classe d'op\'erateurs pseudo-diff\'erentiels sur $\mathbb R^n$. 
\newblock  Math. Nachr. 140 (1989), 127--186.


\bibitem {Pe} J. Peetre.
\newblock A generalization of Courant nodal theorem.
\newblock Math. Scandinavica 5 (1957), 15-20.

\bibitem {Pl} {\AA}.~Pleijel.
\newblock Remarks on Courant's nodal theorem.
\newblock Comm. Pure. Appl. Math. 9 (1956), 543-550.

\bibitem {Pol} I. Polterovich.
\newblock Pleijel's nodal domain theorem for free membranes.
\newblock Proc. Amer. Math. Soc. 137(3) (2009), 1021-1024.

\bibitem{ReSi} M. Reed, B. Simon.
\newblock Method of Modern Mathematical Physics IV: Analysis of operators.
\newblock Academic Press (1978).

\bibitem{Rob} D. Robert.
\newblock Propri\'et\'es spectrales d'op\'erateurs pseudo-diff\'erentiels.
\newblock Comm. in PDE 3(9)  (1978), 755-826.

\bibitem{Roz} G. V. Rosenbljum. 
\newblock The distribution of the discrete spectrum for singular 
differential operators,
\newblock Sov. Math. Dokl. 13 (1972), 245-249.


\bibitem{Sh} 
M.~A. Shubin.
\newblock Asymptotic behaviour of the spectral function.
\newblock In {\em Pseudodifferential Operators and Spectral Theory}, 
  133--173. Springer Berlin Heidelberg, 2001.

 \bibitem{TS}   M. A. Shubin and V. N. Tulovskii.
\newblock 
On the asymptotic distribution of the eigenvalues of pseudodifferential operators in $\mathbb R^n$.
\newblock Mat.  Sb.  92 (134) (1973), 571--588 (in
Russian).

\bibitem{St} S. Steinerberger.
\newblock A geometric uncertainty principle with an application to Pleijel's estimate, 
\newblock Ann. Henri Poincar\'e 15 (2014), no. 12, 2299--2319.


\bibitem{Ta0} H. Tamura.
\newblock Asymptotic formulas with remainder estimates for eigenvalues of Schr\"odinger operators.
\newblock Comm. Partial Differential Equations 7(1), 1-53 (1982).


\bibitem{Ta1} H. Tamura.
\newblock Asymptotic formulas with sharp remainder estimates for bound states of Schr\"odinger I.
\newblock J. Analyse Math. 40 (1981), 166--182.

 \bibitem{Ta2} H. Tamura.
\newblock Asymptotic formulas with sharp remainder estimates for bound states of Schr\"odinger II.
\newblock J. Analyse Math.  41 (1982), 85--108. 


\bibitem {W} H. Weyl.
\newblock \"Uber die asymptotische Verteilung der Eigenwerte.
\newblock Nachrichten der K\"oniglichen Gesellschaft der
Wissenschaften zu G\"ottingen (1911), 110--117.


\end{thebibliography}
\end{document}